\documentclass{article}

\usepackage{amsmath,amssymb,amsfonts,amsthm}

\makeatletter \@addtoreset{equation}{section}
\@addtoreset{theequation}{thesection} \makeatother

\newtheorem{theorem}{Theorem}[section]
\newtheorem{lemma}{Lemma}[section]
\newtheorem{corollary}{Corollary}[section]

\theoremstyle{remark}
\newtheorem{remark}{Remark}[section]

\DeclareMathOperator{\dist}{dist} \DeclareMathOperator{\pos}{pos}

\title{The Hopf boundary point lemma for vector bundle sections}

\author{Artem Pulemotov}

\begin{document}

\maketitle

\begin{center}
Department of Mathematics, Cornell University, \\ 310 Malott Hall,
Ithaca, NY 14853-4201, USA \\ E-mail: artem@math.cornell.edu
\end{center}

\begin{abstract}
The paper establishes a version of the Hopf boundary point lemma
for sections of a vector bundle over a manifold with boundary.
This result may be viewed as a counterpart to the tensor maximum
principle obtained by R.~Hamilton in~1986. Potential applications
include the study of various geometric flows and the construction
of invariant sets for geometric boundary value problems.
\end{abstract}

\section{Introduction}

The present paper concerns the solutions of a second-order partial
differential equation in a vector bundle over a manifold with
boundary. Let us describe the rootage of the considered problem.

The maximum principle for sections of a general vector bundle over a
closed manifold was originally obtained in~\cite{H86}. This
statement is also known as the tensor maximum principle. It proved
to be a powerful implement in the study of the Ricci flow;
see~\cite{CK04}. In particular, it was used to establish important
facts about four-manifolds with nonnegative curvature operator.
Other applications were considered, as well; see, for
instance,~\cite{C92}, \cite{H93}, and~\cite{CN05}. In particular,
the paper~\cite{H93} establishes the important matrix Li-Yau
inequality for solutions of the heat equation.

A specific version of the maximum principle for sections appeared
in~\cite{H82}. This version only applied to 2-tensors. Several
generalizations of the maximum principle for sections can be found
in~\cite{CL04}. We refer to~\cite{C99} and Chapter~4 of~\cite{CK04}
for an overview of relevant results. Once again, we emphasize that
the theory discussed above has been developed on closed manifolds.

The maximum principle for sections may be regarded as an evolution
of the maximum principle for systems of scalar parabolic equations
obtained in~\cite{W75}. It must be noted that the statement
in~\cite{W75} has become a powerful implement in the study of
parabolic systems. In particular, it was applied to the
investigation of the existence and the asymptotic behavior of
solutions. We refer to Chapter~14 of~\cite{Sm94} for several
relevant results and a vast bibliography; some of the references
not mentioned there are~\cite{RW78}, \cite{CS87}, \cite{K00}, and
\cite{AR03}.

An important comment should be made at this point. The maximum
principles discussed above rely on the concept of an invariant
set. The definition of an invariant set for a system of scalar
parabolic equations can be found, for example, in Chapter~14
of~\cite{Sm94}. This definition generalizes easily to cover the
case of an equation for vector bundle sections. We remark that
invariant sets should be viewed from a slightly different
standpoint when the boundary conditions are specified for the
solutions; see, for instance,~\cite{RW78}, \cite{K80},
and~\cite{K00}.

The paper~\cite{Sh96}, being devoted to the study of the Ricci
flow on manifolds with boundary, offers a specific version of the
Hopf boundary point lemma. This version applies to 2-tensors over
a manifold with boundary. In essence, it is an analogue of the
maximum principle for 2-tensors proved in~\cite{H82}. At the same
time, in spite of the fact that the universal maximum principle
for sections suggested in~\cite{H86} is a recognized powerful
tool, no counterparts of this statement have yet been obtained in
the presence of a boundary. Section~\ref{sec result} of the
present paper establishes a general version of the Hopf boundary
point lemma. Our statement applies to sections of a general vector
bundle over a manifold with boundary. It appears to constitute a
comprehensive counterpart to the maximum principle suggested
in~\cite{H86}.

After proving our Hopf lemma for sections, we state three of its
immediate corollaries. They are all closely related to the concept
of an invariant set. The first corollary may be viewed as the basic
maximum principle for sections of a vector bundle over a manifold
with boundary. The second corollary shows that the maximum principle
of~\cite{H86} holds in the presence of a boundary provided that
Neumann-type boundary conditions are imposed on the sections in
question. Such a result is expected to prove useful in the study of
the Ricci flow; cf.~\cite{Sh96}. The third corollary provides an
explicit connection between invariant sets of an equation for vector
bundle sections and the boundary conditions specified for the
solutions. In certain situations, it allows to find an invariant set
for a given boundary value problem. (In one form or another, this
task was addressed in many works; see, for instance,~\cite{K80},
Chapter~14 of~\cite{Sm94}, \cite{K00}, and~\cite{AR03}.)
Alternatively, the third corollary may be used to construct a
relatively sophisticated boundary value problem with a given
invariant set.

Section~\ref{expl for domains} of the present paper explains how
our Hopf lemma for sections applies to a system of scalar
parabolic equations similar to the one studied in~\cite{W75}.
Roughly speaking, we reformulate our statement for sections of a
trivial bundle equipped with the standard connection. An
analogous, although not exactly the same, result for parabolic
systems can be found in~\cite{RW78}. We should note that the
specific framework of Section~\ref{expl for domains} enables us to
refine the smoothness assumptions imposed in Section~\ref{sec
result}.

\section{The Hopf lemma for sections}\label{sec result}

Consider a smooth, compact, connected, oriented Riemannian manifold
$M$ with (possibly empty) boundary $\partial M$. We use the notation
$\nu(x)$ for the outward unit normal to $\partial M$ at the point
$x\in\partial M$. Let $V$ be a vector bundle over $M$. The fiber of
$V$ over $x\in M$ will be denoted by $V_x$. The designation $\pi(v)$
refers to the projection of $v\in V$ onto $M$. We suppose $V$ is
equipped with a fiber metric~$\langle\cdot,\cdot\rangle$.
Let~$\|\cdot\|$ stand for the corresponding norm.

Consider a time-dependent section $f(t,x)$ of the vector bundle
$V$. In what follows, the time parameter $t$ varies through the
interval $[0,T]$ with a fixed $T>0$. Choose a connection $A$ in
$V$ compatible with~$\langle\cdot,\cdot\rangle$. Let $\nabla
f(t,x)$ denote the covariant derivative of $f(t,x)$ with respect
to $A$. We write $\nabla_\chi f(t,x)$ to indicate the application
of $\nabla f(t,x)$ to the tangent vector $\chi\in T_xM$. Employing
the connection $A$ in $V$ and the Levi-Civita connection in the
cotangent bundle $T^*M$, one can define the second covariant
derivative $\nabla^2f(t,x)$. We write
$\nabla^2_{\chi_1,\chi_2}f(t,x)$ to indicate the application of
$\nabla^2f(t,x)$ to the vectors $\chi_1,\chi_2\in T_xM$. The
Laplacian $\Delta$ acts on the section $f(t,x)$ by taking the
trace of $\nabla^2f(t,x)$. We refer to Chapter~4 of~\cite{CK04}
for the details of defining the Laplacian.\footnote{The
denotations $\hat\nabla(t)(\bar\nabla(t)\varphi)(t,x)$ and
$\hat\Delta$ are used in Chapter~4 of~\cite{CK04} for the objects
denoted by $\nabla^2f(t,x)$ and $\Delta$ in the present paper.}

Let $\phi(t,v)$ be a time-dependent mapping of $V$ into itself
such that $\phi(t,v)\in V_{\pi(v)}$ for any $(t,v)\in[0,T]\times
V$. Suppose every compact set $U\subset V$ admits a constant
$C_\phi(U)>0$ satisfying
\begin{align}\label{constant C phi} \|\phi(t,v_1)-\phi(t,v_2)\|\le
C_\phi(U)\|v_1-v_2\|.
\end{align}
The estimate must hold for any $t\in(0,T)$, and any $v_1,v_2\in U$
subject to $\pi(v_1)=\pi(v_2)$. Let $\zeta(t,x)$ be a
time-dependent vector field on $M$. Suppose $f(t,x)$ solves the
second-order equation
\begin{align}\label{heat eq bund}
\frac{\partial}{\partial t}f(t,x)=\Delta
f(t,x)+\nabla_{\zeta(t,x)} f(t,x)+\phi(t,f(t,x))
\end{align}
on $(0,T)\times M$. In particular, $f(t,x)$ must be continuous in
$t\in[0,T]$ and $C^1$-differentiable in $t\in(0,T)$.

Consider a nonempty set $W\subset V$. We assume $W$ is invariant
under the parallel translation with respect to the connection $A$
fixed in $V$. The set $W_x=W\cap V_x$ must be closed and convex in
the fiber $V_x$ for every $x\in M$. When writing $\partial W_x$,
we refer to the boundary of $W_x$ in $V_x$. It should be noted
that $\partial W_x$ is not required to be smooth for any $x\in M$.
Given a point $\omega\in W$ subject to $\omega\in\partial
W_{\pi(\omega)}$, we call $\lambda\in V_{\pi(\omega)}$ a
\emph{supporting vector} for $W$ at $\omega$ if $\|\lambda\|=1$
and the inequality $\langle\lambda,\sigma\rangle \le
\langle\lambda,\omega\rangle$ holds for all $\sigma\in
W_{\pi(\omega)}$. The set of all the supporting vectors for $W$ at
$\omega$ will be denoted by $S_\omega W$. In a sense, the elements
of $S_\omega W$ are outward unit normals to $\partial
W_{\pi(\omega)}$ at $\omega$.

Introduce the notation
\begin{align*}
\dist_Wv=\inf_{\omega\in W_{\pi(v)}}\|v-\omega\|
\end{align*} for $v\in V$. Let $\omega(v)$ be the unique
point in $W_{\pi(v)}$ such that $\dist_Wv=\|v-\omega(v)\|$.
Obviously, $\dist_Wv$ represents the distance between $v\in V$ and
$W_{\pi(v)}$, while $\omega(v)$ is the unique point in
$W_{\pi(v)}$ closest to $v$. We call $(t,x)\in[0,T]\times M$ a
\emph{maximal distance pair} if
\begin{align*} \dist_Wf(t,x)=\sup_{y\in M}\dist_Wf(t,y)>0.
\end{align*} Let $\lambda(v)$ denote the difference $v-\omega(v)$ for $v\in V$.

We are now ready to formulate our Hopf lemma for sections. It
should be remarked that the assumption on the mapping $\phi(t,v)$
in our statement is quite standard. Roughly speaking, we demand
that $\phi(t,v)$ point into $W$ when $v$ is subject to
$v\in\partial W_{\pi(v)}$. This is equivalent to the ``ordinary
differential equation assumption" employed in~\cite{H86}; see
Lemma~4.1 in~\cite{H86}.

\begin{theorem}\label{Str max vbund}
Suppose the solution $f(t,x)$ of equation~\eqref{heat eq bund} and
the mapping $\phi(t,v)$ appearing in the right-hand side of
equation~\eqref{heat eq bund} meet the following requirements:
\begin{enumerate}
\item
The initial value $f(0,x)$ lies in $W$ for all $x\in M$.
\item
The estimate $\langle\lambda,\phi(t,\omega)\rangle\le0$ holds for
any $t\in(0,T)$, any $\omega\in W$ subject to $\omega\in\partial
W_{\pi(\omega)}$, and any supporting vector $\lambda\in S_\omega
W$.
\end{enumerate}
If the value $f(t,x)$ lies outside of $W$ for some
$(t,x)\in(0,T]\times M$, then there exists a maximal distance pair
$(t_{\pos},x_{\pos})\in(0,T)\times\partial M$ such that the
formula
\begin{align}\label{statement}
\left\langle\lambda(f(t_{\pos},x_{\pos})),\nabla_{\nu(x_{\pos})}
f(t_{\pos},x_{\pos})\right\rangle&>0\end{align} holds true.
\end{theorem}

Before proving the theorem, we need to make some preliminary
arrangements. Given a real-valued function $\theta(t)$ on $[0,T)$,
define
\begin{align*}\dot\theta^+(t)=\limsup_{h\to0+}\frac{\theta(t+h)-\theta(t)}h
\end{align*}
for $t\in[0,T)$. The following lemma will be required;
cf.~Lemma~3.1 and Corollary~3.3 in~\cite{H86}, or Lemma~7
in~\cite{CL04}.

\begin{lemma}\label{aux_lip}
Suppose $\theta(t)$ is a nonnegative continuous function on $[0,T)$
with $\theta(0)=0$. Suppose also $\theta(t)$ is not identically $0$
on $[0,T)$. Given a constant $C>0$, there exists a point
$t_C\in(0,T)$ such that $\dot\theta^+(t_C)>C\theta(t_C)$ and
$\theta(t_C)>0$.
\end{lemma}
\begin{proof}
Assume the existence of $C>0$ satisfying the estimate
$\dot\theta^+(t)\le C\theta(t)$ whenever $\theta(t)>0$. Introduce a
new nonnegative continuous function $\eta(t)=e^{-Ct}\theta(t)$.
Clearly, the equality $\eta(0)=0$ holds, and $\dot\eta^+(t)\le0$
whenever $\eta(t)>0$.

Fix $\epsilon_1,\epsilon_2>0$. We will now prove that
$\eta(t)\le\epsilon_1t+\epsilon_2$ for all $t\in[0,T)$. Let $a$ be
the largest possible number in $(0,T]$ such that the inequality
$\eta(t)\le\epsilon_1t+\epsilon_2$ holds on $[0,a)$. (Since
$\eta(0)=0<\epsilon_2$, the set of such numbers is not empty, and
$a$ is well defined.) We claim that $a=T$. Indeed, if $a<T$, then
$\eta(a)=\epsilon_1a+\epsilon_2>0$ by continuity and
\begin{align*}\limsup_{h\to0+}\frac{\eta(a+h)-\eta(a)}h\le0.
\end{align*}
But this implies $\eta(t)\le\epsilon_1t+\epsilon_2$ on
$[0,a+\delta)$ for some $\delta>0$, which contradicts the
definition of $a$.

Thus $\eta(t)\le\epsilon_1t+\epsilon_2$ for all $t\in[0,T)$. Since
this inequality holds for any $\epsilon_1,\epsilon_2>0$, we can
conclude that $\eta(t)$ is identically 0. Hence $\theta(t)$ is
identically 0, which contradicts the suppositions of the lemma.
\end{proof}

\begin{proof}[Proof of Theorem~\ref{Str max vbund}.]

It suffices to carry out the proof assuming $W$ is compact. In order
to justify this statement, fix a number $R>0$ large enough to ensure
that $\|f(t,x)\|<R$ and $\|\omega(f(t,x))\|<R$ for any
$(t,x)\in[0,T]\times M$. Introduce the set $\hat W=\{w\in
W|\,\|w\|\le R\}$. One can verify that $\hat W$ is compact. Clearly,
it is invariant under the parallel translation with respect to $A$,
and its intersection with the fiber $V_x$ is closed and convex in
$V_x$ for every $x\in M$. Let $\kappa(v)$ be a smooth function
acting from $V$ to the interval $[0,1]$. We choose $\kappa(v)$
demanding that $\kappa(v)=1$ when $\|v\|\le R$ and $\kappa(v)=0$
when $\|v\|\ge2R$. Define the time-dependent mapping $\hat\phi(t,v)$
of $V$ into itself by the formula
$\hat\phi(t,v)=\kappa(v)\phi(t,v)$. Estimate~\eqref{aux_lip} is
obviously satisfied for $\hat\phi(t,v)$ with the constant
$C_{\hat\phi}(U)=C_\phi(U)$ when the compact set $U$ is equal to
$f([0,T]\times M)\cup\hat W$. (We note that the proof of the theorem
will not require estimate~\eqref{aux_lip} to hold when $U$ is other
than $f([0,T]\times M)\cup W$.) The section $f(t,x)$ would remain a
solution of equation~\eqref{heat eq bund} if the mapping
$\hat\phi(t,v)$ appeared in the right-hand side of this equation
instead of the mapping $\phi(t,v)$. A straightforward argument
demonstrates that it suffices to prove the theorem with $W$ and
$\phi(t,v)$ replaced by $\hat W$ and $\hat\phi(t,v)$. Therefore,
supposing $W$ is compact does not lead to a loss of generality.

Introduce the function
\begin{align*}
s(t)=\sup_{x\in M}\dist_Wf(t,x)
\end{align*}
for $t\in[0,T]$. Evidently, it is nonnegative. One can show that
$s(t)$ is continuous. Our requirement~1 implies that $s(0)=0$. If
$f(t,x)$ lies outside of $W$ for some $(t,x)\in(0,T]\times M$, then
$s(t)$ is not identically 0 on $[0,T]$. Assuming the assertion of
the theorem fails to hold, we will prove the estimate $\dot
s^+(t)\le Cs(t)$ for a fixed constant $C>0$ and an arbitrary
$t\in(0,T)$ such that $s(t)>0$. Lemma~\ref{aux_lip} would then
provide a contradiction.

Fix a point $t\in(0,T)$ satisfying $s(t)>0$. When $x\in M$ is
subject to $\dist_Wf(t,x)>0$, the equality
\begin{align*}
\dist_Wf(t,x)=\sup_{\omega\in\partial W_x}\sup_{\lambda\in
S_\omega W}\langle\lambda,f(t,x)-\omega\rangle
\end{align*}
holds true. This implies
\begin{align*}
s(t)&=\sup_{(\omega,\lambda)\in\Omega}\langle\lambda,f(t,\pi(\omega))-\omega\rangle,
\\ \Omega&=\left\{(\omega,\lambda)\in V\times V\,\big|\,\omega\in\partial W_{\pi(\omega)},\lambda\in S_\omega
W\right\}.
\end{align*}
The set $\Omega$ is compact in $V\times V$. Therefore, we can
apply Lemma~9 in~\cite{CL04}, see also Lemma~3.5 in~\cite{H86}, to
conclude
\begin{align*}
\dot s^+(t)
&\le\sup_{(\omega,\lambda)\in\Omega'}\frac\partial{\partial r}
\langle\lambda,f(r,{\pi(\omega)})-\omega\rangle|_{r=t},
\\
\Omega'&=\{(\omega,\lambda)\in\Omega\,|\,s(t)=\langle\lambda,f(t,\pi(\omega))-\omega\rangle\}.
\end{align*}

Fix a pair $(\omega,\lambda)\in\Omega'$. For brevity, we write $x$
instead of $\pi(\omega)$. The point $x\in M$ is thus fixed from
now on. Assuming the assertion of the theorem fails to hold, we
will show that $\frac\partial{\partial r}
\langle\lambda,f(r,x)-\omega\rangle|_{r=t}\le Cs(t)$ for a
constant $C>0$ independent of $t$. This would yield the desired
estimate $\dot s^+(t)\le Cs(t)$.

Equation~\eqref{heat eq bund} yields
\begin{align}\label{heat eq functl}
\frac\partial{\partial r}
\langle\lambda,f(r,x)&-\omega\rangle|_{r=t} \notag
\\ &=\bigl\langle\lambda,\Delta
f(t,x)\bigr\rangle+\left\langle\lambda,\nabla_{\zeta(t,x)}
f(t,x)\right\rangle+\bigl\langle\lambda,\phi(t,f(t,x))\bigr\rangle\,.
\end{align} The inclusion $(\omega,\lambda)\in\Omega'$ implies that
$(t,x)$ is a maximal distance pair and the vector $\lambda$
coincides with $\frac{\lambda(f(t,x))}{\|\lambda(f(t,x))\|}$. If the
assertion of the theorem were incorrect, then either $x$ would be in
the interior of $M$ or $\left\langle\lambda,\nabla_{\nu(x)}
f(t,x)\right\rangle$ would be non-positive. Assuming this
alternative, we will estimate each of the three terms in the
right-hand side of equation~\eqref{heat eq functl}.

Let us establish the equality $\left\langle\lambda,\nabla_\chi
f(t,x)\right\rangle=0$ for an arbitrary $\chi\in T_xM$. Obviously,
it would imply
\begin{align}\label{est 1st
term}\left\langle\lambda,\nabla_{\zeta(t,x)}
f(t,x)\right\rangle=0.\end{align} At the first step, we consider a
vector $\chi\in T_x M$ admitting a geodesic segment
$\gamma_\chi(u)$ defined for $u\in[0,\epsilon_\chi]$ in such a way
that $\gamma_\chi(0)=x$ and $\frac{d\gamma_\chi}{du}(0)=\chi$. The
number $\epsilon_\chi$ should be chosen small enough to ensure the
geodesic segment's not intersecting itself. The initial goal is to
show that $\left\langle\lambda,\nabla_{\chi}
f(t,x)\right\rangle\le0$.

For the sake of brevity, we write $\gamma(u)$ instead
of~$\gamma_\chi(u)$ and $\epsilon$ instead of~$\epsilon_\chi$. One
can extend the vectors $\lambda$ and $\omega$ to parallel~(with
respect to the connection~$A$) sections $\lambda'(\gamma(u))$ and
$\omega'(\gamma(u))$ of the bundle $V$ defined along $\gamma(u)$.
The covariant derivatives of $\lambda'(\gamma(u))$ and
$\omega'(\gamma(u))$ with respect to $A$ at the point
$x=\gamma(0)$ exist in the direction of $\chi$. Writing
$\nabla_{\chi}\lambda'(x)$ and $\nabla_{\chi}\omega'(x)$ for these
covariant derivatives, we can easily see that
$\nabla_{\chi}\lambda'(x)=0$ and $\nabla_{\chi}\omega'(x)=0$.

Introduce the function
$g(u)=\left\langle\lambda'(\gamma(u)),f(t,\gamma(u))-\omega'(\gamma(u))\right\rangle$
on $[0,\epsilon]$. Obviously, $g(0)=s(t)$. Using the fact that the
parallel transport is an isometry of the fibers, one proves
$\omega'(\gamma(u))\in\partial W_{\gamma(u)}$ and
$\lambda'(\gamma(u))\in S_{\omega'(\gamma(u))}W$ for any
$u\in[0,\epsilon]$. These inclusions imply the inequality
\begin{align*}
g(0)=s(t)\ge\langle\lambda'(\gamma(u)),f(t,\gamma(u))-\omega'(\gamma(u))\rangle=g(u)
\end{align*}
for any $u\in[0,\epsilon]$. As a consequence, the function $g(u)$
has a maximum at $0$, and the one-sided derivative
$\frac{dg}{du}(0)$ is non-positive. Since the connection $A$ is
compatible with the fiber metric, we have the formula
\begin{align*}
\left\langle\lambda,\nabla_{\chi}
f(t,x)\right\rangle&=\left\langle\nabla_{\chi}\lambda'(x),f(t,x)\right\rangle+\left\langle\lambda,\nabla_{\chi}
f(t,x)\right\rangle \\ &=\frac\partial{\partial
u}\left.\left\langle\lambda'(\gamma(u)),f(t,\gamma(u))\right\rangle\right|_{u=0}=\frac{dg}{du}(0).
\end{align*}
Hence $\left\langle\lambda,\nabla_{\chi}f(t,x)\right\rangle\le0$.

Choose an orthonormal basis $\{e_1,\ldots,e_n\}$ of the tangent
space $T_xM$. We will show that $\left\langle\lambda,\nabla_{e_k}
f(t,x)\right\rangle=0$ for any $k=1,\ldots,n$. Suppose $x$ lies in
the interior of $M$. Then a geodesic segment $\gamma_{e_k}(u)$, the
parameter $u$ varying through~$[0,\epsilon_{e_k}]$, subject to
$\gamma_{e_k}(0)=x$ and $\frac{d\gamma_{e_k}}{du}(0)=e_k$ exists for
any $k=1,\ldots,n$. As a consequence, the scalar products
$\left\langle\lambda,\nabla_{e_k}f(t,x)\right\rangle$ are
non-positive. Substituting $-e_k$ for $e_k$ and repeating the
argument, we conclude that the scalar products
$\left\langle\lambda,\nabla_{e_k}f(t,x)\right\rangle$ are also
nonnegative. Thus
$\left\langle\lambda,\nabla_{e_k}f(t,x)\right\rangle=0$ for any
$k=1,\ldots,n$.

Suppose $x$ lies in the boundary of $M$. Without loss of generality,
we assume $e_n$ coincides with the inward normal to the boundary of
$M$. It is easy to verify the existence of a geodesic segment
$\gamma_{e_n}(u)$ defined for $u\in[0,\epsilon_{e_n}]$ in such a way
that $\gamma_{e_n}(0)=x$ and $\frac{d\gamma_{e_n}}{du}(0)=e_n$.
Consequently, the scalar product $\left\langle\lambda,\nabla_{e_n}
f(t,x)\right\rangle$ is non-positive. At the same time, our
hypothesis implies that $\bigl\langle\lambda,\nabla_{e_n}
f(t,x)\bigr\rangle=-\left\langle\lambda,\nabla_{\nu(x)}
f(t,x)\right\rangle$ is nonnegative. Thus
$\left\langle\lambda,\nabla_{e_n}f(t,x)\right\rangle=0$. Provided
$n\ge2$, we now prove that
$\left\langle\lambda,\nabla_{e_k}f(t,x)\right\rangle=0$ for
$k=1,\ldots,n-1$. The situation is slightly more complicated here
because a geodesic emanating from $x$ in the direction of $e_k$ does
not necessarily exist. In order to overcome this problem, we will
carry out an approximation procedure. Namely, fix a sequence
$(e_k^m)_{m=1}^\infty$ converging to $e_k$ for every
$k=1,\ldots,n-1$. We choose the vectors $e_k^m$ demanding that the
scalar product of $e_k^m$ and $e_n$ with respect to the Riemannian
metric in $M$ be strictly positive. Given $k$ and $m$, it is easy to
verify the existence of a geodesic segment $\gamma_{e_k^m}(u)$, the
parameter $u$ varying through $\left[0,\epsilon_{e_k^m}\right]$,
subject to $\gamma_{e_k^m}(0)=x$ and
$\frac{d\gamma_{e_k^m}}{du}(0)=e_k^m$. As a consequence,
$\left\langle\lambda,\nabla_{e_k^m} f(t,x)\right\rangle\le0$. The
convergence of $(e_k^m)_{m=1}^\infty$ to $e_k$ then implies
$\left\langle\lambda,\nabla_{e_k} f(t,x)\right\rangle\le0$.
Substituting $-e_k$ for $e_k$ and repeating the argument, we
conclude that
$\left\langle\lambda,\nabla_{e_k}f(t,x)\right\rangle\ge0$. Thus
$\left\langle\lambda,\nabla_{e_k}f(t,x)\right\rangle=0$ for
$k=1,\ldots,n-1$.

By virtue of the established equalities,
$\left\langle\lambda,\nabla_{\chi}f(t,x)\right\rangle=0$ for an
arbitrary $\chi\in T_xM$. This clearly proves formula~\eqref{est
1st term}.

Our next goal is to obtain the estimate
\begin{align}\label{est Laplacian}\left\langle\lambda,\Delta
f(t,x)\right\rangle\le0.\end{align} As before, consider a vector
$\chi\in T_x M$ admitting a geodesic segment $\gamma_\chi(u)$
defined for $u\in[0,\epsilon_\chi]$ in such a way that
$\gamma_\chi(0)=x$ and $\frac{d\gamma_\chi}{du}(0)=\chi$. The
number $\epsilon_\chi$ should be small enough to ensure the
absence of self-intersections. We now show that
$\left\langle\lambda,\nabla^2_{\chi,\chi}f(t,x)\right\rangle\le0$.
This would provide us with a basis for the proof of
estimate~\eqref{est Laplacian}.

Again, we write $\gamma(u)$ instead of $\gamma_\chi(u)$ and
$\epsilon$ instead of $\epsilon_\chi$. It will be convenient to
use the notation $\gamma'(u)$ for $\frac{d\gamma}{du}(u)$. A
parallel section $\lambda'(\gamma(u))$ of the bundle $V$ along
$\gamma(u)$ has been introduced above. The covariant derivative of
this section with respect to the connection $A$ at the point
$\gamma(u)$ exists in the direction of $\gamma'(u)$ for any
$u\in[0,\epsilon)$. Writing
$\nabla_{\gamma'(u)}\lambda'(\gamma(u))$ for this covariant
derivative, we can easily see that
$\nabla_{\gamma'(u)}\lambda'(\gamma(u))=0$ for any
$u\in[0,\epsilon)$.

Since $A$ is compatible with the fiber metric, the equality
\begin{align*}\left\langle\lambda,\nabla_{\chi,\chi}^2f(t,x)\right\rangle&=
\bigl\langle\nabla_\chi\lambda'(x),\nabla_\chi f(t,x)\bigr\rangle
+\left\langle\lambda,\nabla^2_{\chi,\chi}f(t,x)\right\rangle
\\ &=\frac\partial{\partial
u}\left\langle\lambda'(\gamma(u)), \nabla_{\gamma'(u)}
f(t,\gamma(u))\right\rangle\big|_{u=0}
\\ &=\frac\partial{\partial
u}\bigl(\left\langle\nabla_{\gamma'(u)}\lambda'(\gamma(u)),f(t,\gamma(u))\right\rangle
\\ &\hphantom{=}~+\left\langle\lambda'(\gamma(u)),\nabla_{\gamma'(u)}
f(t,\gamma(u))\right\rangle\bigr)\big|_{u=0} \\
&=\frac{\partial^2}{\partial
u^2}\left.\left\langle\lambda'(\gamma(u)),f(t,\gamma(u))\right\rangle\right|_{u=0}
\\ &=\frac{d^2g}{du^2}(0)\end{align*} holds true. The introduced above function $g(u)$ has a maximum at 0. It has
been proven that $\frac{dg}{du}(0)=\left\langle\lambda,\nabla_\chi
f(t,x)\right\rangle=0$. Hence $\frac{d^2g}{du^2}(0)\le0$, which
yields
$\left\langle\lambda,\nabla^2_{\chi,\chi}f(t,x)\right\rangle\le0$.

Suppose $x$ lies in the interior of $M$. Then every vector from
the chosen above basis $\{e_1,\ldots,e_n\}$ appears as a tangent
vector for a certain geodesic segment emanating from $x$. As a
consequence,
$\left\langle\lambda,\nabla^2_{e_k,e_k}f(t,x)\right\rangle\le0$
for every $k=1,\ldots,n$.

Suppose $x$ lies in the boundary of $M$. Recall that $e_n$ is
assumed to coincide with the inward normal to the boundary of $M$.
As mentioned before, $e_n$ appears as a tangent vector for a
certain geodesic segment emanating from $x$. Therefore,
$\left\langle\lambda,\nabla^2_{e_n,e_n}f(t,x)\right\rangle\le0$.
Provided $n\ge2$, we can approximate the other basis vectors with
the previously fixed sequences $(e_k^m)_{m=1}^\infty$ to conclude
that
$\left\langle\lambda,\nabla^2_{e_k,e_k}f(t,x)\right\rangle\le0$
for every $k=1,\ldots,n-1$.

According to the definition of the Laplacian,
\begin{align*}\langle\lambda,\Delta
f(t,x)\rangle=\sum_{k=1}^n\left\langle\lambda,\nabla^2_{e_k,e_k}f(t,x)\right\rangle.\end{align*}
By virtue of the established inequalities, all the terms in the
right-hand side are non-positive. This clearly implies
formula~\eqref{est Laplacian}.

Finally, let us prove the estimate
\begin{align}\label{3rd term est}
\langle\lambda,\phi(t,f(t,x))\rangle\le Cs(t)
\end{align}
with a constant $C>0$ independent of $t$. The vector $\lambda$
belongs to $S_\omega W$. It must also belong to
$S_{\omega(f(t,x))}W$, although $\omega$ does not necessarily
coincide with $\omega(f(t,x))$. (Recall that $\omega(f(t,x))$
stands for the unique point in $W_x$ closest to $f(t,x)$.) In
accordance with our requirement~2,
$\langle\lambda,\phi(t,\omega(f(t,x)))\rangle \le0$. Hence the
estimate
\begin{align*}
\langle\lambda,\phi(t,f(t,x))\rangle&\le
\langle\lambda,\phi(t,f(t,x))\rangle-\langle\lambda,\phi(t,\omega(f(t,x)))\rangle
\\ &\le\|\phi(t,f(t,x))-\phi(t,\omega(f(t,x)))\| \\ &\le C\|f(t,x)-\omega(f(t,x))\|=Cs(t)
\end{align*}
holds with the constant $C>0$ equal to the constant
$C_\phi(f([0,T]\times M)\cup W)>0$ given by
formula~\eqref{constant C phi}. This concludes the proof
of~\eqref{3rd term est}. Remark that the argument we used does not
depend on whether $x$ is in the boundary of $M$ or in the interior
of $M$.

Equation~\eqref{heat eq functl} now provides
$\frac\partial{\partial r}
\langle\lambda,f(r,x)-\omega\rangle|_{r=t}\le Cs(t)$. As mentioned
before, this inequality implies $\dot s^+(t)\le Cs(t)$, which is
impossible in view of Lemma~\ref{aux_lip}.
\end{proof}

\begin{remark}
The assumption on the mapping $\phi(t,v)$ imposed by the theorem may
be slightly refined. Namely, it suffices to demand that the estimate
$\langle\lambda,\phi(t,\omega)\rangle\le0$ hold when $\omega$ is
equal to $\omega(f(t,x))$ and $\lambda$ is equal to
$\lambda(f(t,x))$ for all the maximal distance pairs
$(t,x)\in(0,T)\times M$.
\end{remark}

\begin{remark}
If the boundary of $M$ is empty, then the suppositions of the
theorem cannot be satisfied simultaneously. In this case,
requirements~1 and~2 ensure that $f(t,x)$ cannot lie outside of
$W$. This fact is essentially equivalent to the maximum principle
obtained in~\cite{H86}.
\end{remark}

\begin{remark}
The theorem would prevail if the Riemannian metric in $M$ and the
connection $A$ fixed in $V$ depended on the time parameter
$t\in[0,T]$. Of course, then we would have to modify some of the
assumptions imposed above. Firstly, the connection $A(t)$ fixed in
$V$ at time $t$ would be required to be compatible with the fiber
metric~$\langle\cdot,\cdot\rangle$ for all $t\in(0,T)$. Secondly,
the set $W$ would have to be invariant under the parallel
translation with respect to $A(t)$ for all $t\in(0,T)$. The
details of defining the Laplacian and writing down
equation~\eqref{heat eq bund} in the situation under discussion
can be found in Chapter~4 of~\cite{CK04}. The covariant derivative
and the outward normal in formula~\eqref{statement} would have to
be computed with respect to the connection $A(t_{\pos})$ and the
Riemannian metric in $M$ at time~$t_{\pos}$.
\end{remark}

We will now formulate three immediate corollaries of
Theorem~\ref{Str max vbund}. The following statement may be viewed
as the basic maximum principle for sections of a vector bundle
over a manifold with boundary.

\begin{corollary}\label{Cor Weak MP}
Suppose the solution $f(t,x)$ and the mapping $\phi(t,v)$ meet
requirements~1 and~2 of Theorem~\ref{Str max vbund}. If $f(t,x)$
lies in $W$ for all $(t,x)\in(0,T)\times\partial M$, then $f(t,x)$
lies in $W$ for all $(t,x)\in[0,T]\times M$.
\end{corollary}

The following statement shows that the maximum principle
of~\cite{H86} holds for $f(t,x)$ provided that Neumann-type boundary
conditions are imposed.

\begin{corollary}\label{Cor Neumann MP}
Suppose the solution $f(t,x)$ and the mapping $\phi(t,v)$ meet
requirements~1 and~2 of Theorem~\ref{Str max vbund}. If the
boundary condition
\begin{align*}\nabla_{\nu(x)}f(t,x)=0\end{align*}
is satisfied for all $(t,x)\in(0,T)\times\partial M$, then
$f(t,x)$ lies in $W$ for all $(t,x)\in[0,T]\times M$.
\end{corollary}

Let $\bar\lambda(v)$ be a mapping of $V$ into itself such that
$\bar\lambda(v)\in V_{\pi(v)}$ for any $v\in V$. The following
statement establishes an explicit connection between invariant
sets of equation~\eqref{heat eq bund} and the boundary conditions
specified for the solutions.

\begin{corollary}\label{cor_inv sets}
Suppose the solution $f(t,x)$ and the mapping $\phi(t,v)$ meet
requirements~1 and~2 of Theorem~\ref{Str max vbund}. Suppose also
$\bar\lambda(v)=\lambda(v)$ for any $v\in V$ lying outside of $W$.
If the boundary condition
\begin{align*}
\left\langle\bar\lambda(f(t,x)),\nabla_{\nu(x)}
f(t,x)\right\rangle=0
\end{align*}
is satisfied for all $(t,x)\in(0,T)\times\partial M$, then
$f(t,x)$ lies in $W$ for all $(t,x)\in[0,T]\times M$.
\end{corollary}

It should be noted that both Corollary~\ref{Cor Weak MP} and
Corollary~\ref{Cor Neumann MP} can be deduced from
Corollary~\ref{cor_inv sets}.

\section{Systems of parabolic equations}\label{expl for domains}

We will now explain how Theorem~\ref{Str max vbund} applies to a
parabolic system similar to the one studied in~\cite{W75}. We
remark that dealing with a parabolic system rather than an
equation for vector bundle sections enables us to refine the
smoothness assumptions imposed in Section~\ref{sec result}.

Let $M$ be the closure of a bounded domain in $\mathbb R^n$ with
$C^1$-differentiable boundary $\partial M$. We use the designation
$\nu(x)$ for the outward unit normal to $\partial M$ at the point
$x\in\partial M$. Differentiation with respect to $x\in M$ in the
direction of $\nu(x)$ will be denoted by
$\frac\partial{\partial\nu}$. Let $x_1,\ldots,x_n$ be the standard
coordinates in $\mathbb R^n$.

Consider a collection of time-dependent real-valued functions
$f_i(t,x)$ on $M$ indexed by $i=1,\ldots,m$. The time parameter
$t$ varies through the interval $[0,T]$. We write $f(t,x)$ for the
vector $(f_1(t,x),\ldots,f_m(t,x))$. Roughly speaking, $f(t,x)$
appears as a time-dependent section of the product bundle
$M\times\mathbb R^m$.

Let $\phi_i(t,x,v)$ be a collection of time-dependent real-valued
functions on $M\times\mathbb R^m$ indexed by $i=1,\ldots,m$.
Again, $\phi(t,x,v)$ stands for the vector
$(\phi_1(t,x,v),\ldots,\phi_m(t,x,v))$. The
denotations~$\langle\cdot,\cdot\rangle$ and~$\|\cdot\|$ refer to
the standard Euclidean scalar product and the standard Euclidean
norm in~$\mathbb R^m$. We demand that every compact set
$U\subset\mathbb R^m$ admit a constant $C_\phi(U)>0$ such that
\begin{align*}\|\phi(t,x,v_1)-\phi(t,x,v_2)\|\le
C_\phi(U)\|v_1-v_2\|.
\end{align*}
The estimate must hold for any $t\in(0,T)$, any $x\in M$, and any
pair of vectors $v_1,v_2\in U$. Roughly speaking, $\phi(t,x,v)$
appears as a time-dependent mapping of the bundle $M\times\mathbb
R^m$ into itself. The above inequality may then be viewed as a
special case of inequality~\eqref{constant C phi}.

Fix time-dependent real-valued functions $\zeta_j(t,x)$ on $M$ for
$j=1,\ldots,n$. We suppose $f_i(t,x)$ solves the second-order
equation
\begin{align}\label{parab system}
\frac\partial{\partial
t}f_i(t,x)=\sum_{j=1}^n\frac{\partial^2}{\partial x_j^2}f_i(t,x)
+\sum_{j=1}^n\zeta_j(t,x)\frac\partial{\partial x_j}f_i(t,x)
&+\phi_i(t,x,f(t,x))
\end{align}
on $(0,T)\times M$ for every $i=1,\ldots,m$. In particular,
$f_i(t,x)$ must be continuous in $t\in[0,T]$, $C^1$-differentiable
in $t\in(0,T)$, and $C^2$-differentiable in $x\in M$. The
collection of equations~\eqref{parab system} for $i=1,\ldots,m$
may be viewed as a special case of equation~\eqref{heat eq bund}.

Consider a nonempty closed convex set $W\subset\mathbb R^m$. Given
a point $\omega\in\partial W$, we call $\lambda\in\mathbb R^m$ a
\emph{supporting vector} for $W$ at $\omega$ if $\|\lambda\|=1$
and the inequality
$\langle\lambda,\sigma\rangle\le\langle\lambda,\omega\rangle$
holds for all $\sigma\in W$. Let $S_\omega W$ signify the set of
all the supporting vectors for $W$ at $\omega$. We write
$\dist_Wv$ to denote the standard Euclidean distance between
$v\in\mathbb R^m$ and $W$. The designation $\omega(v)$ refers to
the unique point in $W$ closest to $v$. We call
$(t,x)\in[0,T]\times M$ a \emph{maximal distance pair} if
\begin{align*} \dist_Wf(t,x)=\sup_{y\in M}\dist_Wf(t,y)>0.
\end{align*}
Let $\lambda(v)$ denote the difference $v-\omega(v)$ for
$v\in\mathbb R^m$.

Following the proof of Theorem~\ref{Str max vbund}, one can obtain
the following result.

\begin{theorem}

Suppose the vector-valued functions $f(t,x)$ and $\phi(t,x,v)$
related by equations~\eqref{parab system} for $i=1,\ldots,m$ meet
the following requirements:
\begin{enumerate}
\item
The initial value $f(0,x)$ lies in $W$ for all $x\in M$.
\item
The estimate $\langle\lambda,\phi(t,x,\omega)\rangle\le0$ holds
for any $(t,x,\omega)\in(0,T)\times M\times\partial W$ and any
supporting vector $\lambda\in S_\omega W$.
\end{enumerate}
If the value $f(t,x)$ lies outside of $W$ for some
$(t,x)\in(0,T]\times M$, then there exists a maximal distance pair
$(t_{\pos},x_{\pos})\in(0,T)\times\partial M$ such that the
formula
\begin{align*}
\frac\partial{\partial\nu}\left\langle\lambda(f(t_{\pos},x_{\pos})),
f(t_{\pos},x)\right\rangle|_{x=x_{\pos}}&>0\end{align*} holds
true.
\end{theorem}

\section*{Acknowledgements}
I express my profound gratitude to Prof.~Leonard Gross for his
support and numerous productive discussions.

\end{document}